# How to Detect Outliers in Data Envelopment Analysis by Kourosh and Arash Method

**Dariush Khezrimotlagh**[*]

Department of Applied Statistics, Faculty of Economics and Administration, University of Malaya, Kuala Lumpur, Malaysia, dariush@um.edu.my, July 01, 2013.

## Abstract

One of the concerns about using non-parametric estimators such as Data Envelopment Analysis (DEA), is the presence of outliers. There are a good number of studies that mention this assessment in the literature of DEA, however, there is no clear definition to identify what outliers are in DEA. Moreover, most of the studies have used additional procedures which have high computational complexities. This paper proposes a suitable definition to identify outliers as well as a simple methodology to illustrate how DEA, by using Kourosh and Arash Method (KAM), is easily able to detect outliers without using additional technologies and their computational complexities. The methodology of detecting outliers by KAM is represented with an example which was used in previous research to depict DEA's weakness of detecting outliers. The results clearly reject this claim that DEA is not able to detect outliers.

*Keywords*: DEA, KAM, Outliers, Ranking, Efficiency.

## 1. Introduction

Data Envelopment Analysis (DEA) is a popular non-parametric technique to measure the relative efficiency of homogenous Decision Making Units (DMUs) with multiple inputs and multiple outputs. It was proposed by Charnes *et al.* (1978), and has dramatically been developed in the last three decades. A good illustration on DEA literature and its conventional models can be seen in Ray (2005), Coelli *et al.* (2005) and Cooper *et al.* (2007).

Seiford and Thrall (1990) clearly illustrated the advantages of DEA to provide new insights and additional information not available from conventional econometric methods. However, there have been some concerns about using DEA technique, such as the influence of presenting outliers in data (Sexton 1986). Detecting outliers, leverage and deviant observations have been one of the controversial objectives in statistics and a good discussion on this issue can be found in Andrews and Pregibon (1977) and Hodge and Austin (2004).

Grosskopf and Valdmanis (1987) and Charnes and Neralić (1990) were the first researchers to attempt rectifying the DEA weakness of detecting outliers. Soon after, Wilson (1993) demonstrated that these methods do not detect inefficient outliers, and proposed a methodology which was also improved for applying large samples by Wilson (1995). Dusansky and Wilson (1995) introduced a methodology for detecting outliers in linear programming efficiency models by deleting n-tuples of observations in order to examine the affected measured efficiency for the remaining observations. Simar (1996) proposed some trails to introduce stochastic noise in DEA models regarding a survey of

---

[*] e-mail address: khezrimotlagh@gmail.com

Grosskopf (1995) and the paper of Banker (1996). Pastor *et al.* (1999) suggested a methodology to detect influential observations which was based on the Andersen and Petersen (1993) method. Simar and Wilson (2000) suggested a general methodology to bootstrap frontier models. Ondrich and Ruggiero (2002) re-sampled the technique of jackknifing to detect outliers, and argued that calculation of the standard deviation is not meaningful for distinguishing different types of outliers. Simar (2003) summarized the results of order-m method proposed by Cazals *et al.* (2002), and showed how it can detect outliers. Chen and Johnson (2006) proposed a model to examine the effect of outliers on the frontier. Banker and Chen (2006) discussed two alternative uses of a super-efficiency model according to the proposed model of Andersen and Petersen (1993) to detect outliers and ranking DMUs. Simar (2007) presented the Hall-Simar procedure inclusive of some simulated examples to detect outliers and the presence of noise. Johnson and McGinnis (2008) considered the efficient and inefficient frontiers to identify outliers based on the works of Wilson (1995), Simar (2003) and Simar and Wilson (2007). Tran *et al.* (2010) introduced a method based on two scalar measures to detect outliers from the set of technically efficient firms. Chen and Johnson (2010) proposed a set of axioms and developed a consistent approach with the axioms for identifying outliers which influence both efficiency estimates and DEA post analysis. Kuosmanen and Johnson (2010) attempted to indicate that DEA can be recast as non-parametric least squares regression to bridge the conceptual gap. Witte and Marques (2010) surveyed several outlier detection methodologies to suggest a procedure to detect outliers from exogenous influences and undesired results of the Portuguese drinking water sector. Emrouznejad and Witte (2010) also discussed this issue in their proposed COOPER-framework to carry non-parametric projects out due to impossibility of assessing the performance entities when the number of factors are increased or when data have to be compiled from several sources. Simar and Zelenyuk (2011) extended the work of Simar (2007) to introduce a methodology for noise in non-parametric frontier models which is based on the theory of local maximum likelihood. Some applications on these issues can also be found in Sengupta (1991), Fried *et al.* (2002), Mahlberg and Raveh (2002), Drake and Simar (2002), Estache *et al.* (2004), de Sousa and Stošić (2005), Dulá (2008), De Jorge Moreno and Sanz- Triguero (2010), Chung (2011) and Horta et al. (2012).

*The above selected studies in the literature of DEA have been done during last three decades due to claim of existing limitations in the conventional DEA models to detect outliers, and unfortunately the claim has been remained, and is still cited in the literature.*

This paper illustrates that such shortcomings is not valid in DEA while Kourosh and Arash Model (KAM) is applied. The paper argues the presence of outliers, and illustrates DEA as a powerful technique and a sufficient tool based on KAM structure to measure the performance evaluation of DMUs with multiple inputs and multiple outputs inclusive real, integer and controllable/non-controllable data. KAM was recently proposed by Khezrimotlagh et al. (2013a) to improve the foundation of DEA and cover many DEA subjects. It is a flexible model which can be adapted for many aims when the weights of factors are available or unknown. KAM significantly ranks both technically efficient and inefficient DMUs, and benchmarks DMUs towards the economical part of the frontier. It is easily solved without the computational complexity of current hybrid DEA methodologies. The properties of KAM can be seen in Khezrimotlagh *et al.* (2012a-f, 2013a-d) and Khezrimotlagh (2014a-b).

The rest of this paper is organized in five sections. Section 2 is a short background on KAM. A discussion on outliers is represented in Section 3 with a proposed definition to identify outliers and classified homogeneous DMUs. The robustness of KAM to detect outliers is exemplified with a simple numerical example in Section 4 according to the example of Simar (2007) and Simar and Zelenyuk (2011). The paper is concluded in the last section. Simulations are also performed using Microsoft Excel Solver, as it required simple linear programming.

## 2. Background on KAM

Suppose that there are $n$ DMUs ($\mathbf{DMU}_i, i = \mathbf{1,2,\ldots,}n$) with $m$ non-negative inputs ($x_{ij}, j = \mathbf{1,2,\ldots,}m$) and $p$ non-negative outputs ($y_{ik}, k = \mathbf{1,2,\ldots,}p$), such that, at least one of the inputs and one of the outputs of each DMU are not zero, and for every $i$ there is a $j$ such that $x_{ij} \neq \mathbf{0}$ and also for every $i$ there is a $k$ such that $y_{ik} \neq \mathbf{0}$. Khezrimotlagh *et al*. (2013a) defined an efficient tape instead of the Farrell frontier, given by:

**Definition 1**: An $\epsilon$-efficient tape is the geometric locus of points in $\mathbb{R}_+^{m+p}$ which their distances of the Farrell frontier do not greater than $\varepsilon_j^-$ of inputs and $\varepsilon_k^+$ of outputs, for $j = \mathbf{1,2,\ldots,}m$ and $k = \mathbf{1,2,\ldots,}p$, where $\epsilon = (\epsilon^-, \epsilon^+) \in \mathbb{R}_+^{m+p}$.

The linear KAM for an appropriate $\epsilon$ while ($\mathbf{DMU}_l, l = \mathbf{1,2,\ldots,}n$) is under evaluation is as follows:

$$\begin{aligned}
&\mathbf{max}\ \sum_{j=1}^{m} w_j^- s_{lj}^- + \sum_{k=1}^{p} w_k^+ s_{lk}^+, &&(1)\\
&\text{Such that}\\
&\sum_{i=1}^{n} \lambda_i x_{ij} + s_{lj}^- = x_{lj} + \varepsilon_j^-,\ \forall j, \quad \sum_{i=1}^{n} \lambda_i y_{ik} - s_{lk}^+ = y_{lk} - \varepsilon_k^+, \forall k,\\
&x_{lj} - s_{lj}^- \geq \mathbf{0},\ \forall j, \qquad\qquad\qquad y_{lk} + s_{lk}^+ - \mathbf{2}\varepsilon_k^+ \geq \mathbf{0}, \forall k,\\
&\sum_{i=1}^{n} \lambda_i = \mathbf{1}, \qquad\qquad\qquad\qquad \lambda_i \geq \mathbf{0},\ \forall i,\\
&s_{lj}^- \geq \mathbf{0},\ \forall j, \qquad\qquad\qquad\qquad s_{lk}^+ \geq \mathbf{0},\ \forall k.
\end{aligned}$$

The targets of KAM, that is, the highest $(\widehat{x^*}, \widehat{y^*})$, the best technical $(x^*, y^*)$ and the lowest $(\widetilde{x^*}, \widetilde{y^*})$ efficient targets with $\epsilon$-Degree of Freedom ($\epsilon$-DF), the highest $\widehat{KA_\epsilon^{*l}}$, the best technical $KA_\epsilon^{*l}$ and the lowest $\widetilde{KA_\epsilon^{*l}}$ efficiency scores as well as sensitivity score $S_\epsilon^{*l}$ with $\epsilon$-DF, respectively are given by:

$$\begin{cases} \widehat{x_{lj}^*} = x_{lj} - s_{lj}^{-*}, \forall j, \\ \widehat{y_{lk}^*} = y_{lk} + s_{lk}^{+*}, \forall k, \end{cases} \qquad \widehat{KA_\epsilon^{*l}} = \frac{\sum_{k=1}^{p} w_k^+ y_{lk} / \sum_{j=1}^{m} w_j^- x_{lj}}{\sum_{k=1}^{p} w_k^+ \widehat{y_{lk}^*} / \sum_{j=1}^{m} w_j^- \widehat{x_{lj}^*}},$$

$$\begin{cases} x_{lj}^* = x_{lj} - s_{lj}^{-*} + \varepsilon_j^-, \forall j, \\ y_{lk}^* = y_{lk} + s_{lk}^{+*} - \varepsilon_k^+, \forall p, \end{cases} \qquad KA_\epsilon^{*l} = \frac{\sum_{k=1}^{p} w_k^+ y_{lk} / \sum_{j=1}^{m} w_j^- x_{lj}}{\sum_{k=1}^{p} w_k^+ y_{lk}^* / \sum_{j=1}^{m} w_j^- x_{lj}^*},$$

$$\begin{cases} \widetilde{x_{lj}^*} = x_{lj} - s_{lj}^{-*} + \mathbf{2}\varepsilon_j^-, \forall j, \\ \widetilde{y_{lk}^*} = y_{lk} + s_{lk}^{+*} - \mathbf{2}\varepsilon_k^+, \forall p, \end{cases} \qquad \widetilde{KA_\epsilon^{*l}} = \frac{\sum_{k=1}^{p} w_k^+ y_{lk} / \sum_{j=1}^{m} w_j^- x_{lj}}{\sum_{k=1}^{p} w_k^+ \widetilde{y_{lk}^*} / \sum_{j=1}^{m} w_j^- \widetilde{x_{lj}^*}},$$

$$S_\epsilon^{*l} = \frac{\widehat{KA_\epsilon^{*l}}}{\widetilde{KA_\epsilon^{*l}}}.$$

After the optimization the following definition is also able to identify the efficient DMUs from the technically efficient ones.

**Definition 2:** A technical efficient DMU is KAM efficient with $\epsilon$-DF in inputs and outputs if $KA_0^* - KA_\varepsilon^* \leq \delta$. Otherwise, it is inefficient with $\epsilon$-DF in inputs and outputs. The proposed amount for $\delta$ is '$\mathbf{10^{-1}\epsilon}$' or '$\epsilon/(m+p)$' (Khezrimotlagh *et al.* 2013a).

The non-linear KAM is given by:

$$KA_\epsilon^{**l} = \min \frac{1 + \sum_{j=1}^{m} W_j^- (E_j^- - s_j^-)}{1 + \sum_{k=1}^{p} W_k^+ (s_k^+ - E_k^+)}, \quad (2)$$

Such that
$\sum_{i=1}^{n} \lambda_i x_{ij} + s_{lj}^- = x_{lj} + E_j^-, \quad \forall j, \quad \sum_{i=1}^{n} \lambda_i y_{ik} - s_{lk}^+ = y_{lk} - E_k^+, \forall k,$
$\sum_{i=1}^{n} \lambda_i = 1; \qquad\qquad\qquad\qquad \lambda_i \geq 0, \quad \forall i,$
$s_{lj}^- \geq 0, \quad \forall j, \qquad\qquad\qquad\qquad s_{lk}^+ \geq 0, \quad \forall k.$

where $W_j^- = w_j^- / \sum_{j=1}^{m} w_j^- x_{lj}$, $W_k^+ = w_k^+ / \sum_{k=1}^{p} w_k^+ y_{lk}$, $E_j^- = \varepsilon_j^- / w_j^-$ and $E_k^+ = \varepsilon_k^+ / w_k^+$. The score of non-linear KAM is denoted by $KA_\epsilon^{**l}$, and it is not greater than the best technical efficiency score of linear KAM denoted by $KA_\epsilon^{*l}$ (Khezrimotlagh *et al.* 2013c). For more information about how KAM is applied see Khezrimotlagh (2014b).

## 3. Outliers

There are several studies and additional methodologies to detect outliers in DEA, however, there is no clear definition to identify outliers in literature of DEA. Wilson (1995) quoted that "outliers are observations that do not fit in with the pattern of the remaining data points, and are not at all typical of the rest of the data (Gunst and Mason, 1980)". He also defined influential observations as "those sample observations which play a relatively large role in determining estimated efficiency scores for at least some other observations in the observed sample". Barnett and Lewis (1984) defined that "outliers are inconsistent observations with the remainder set of data". Davies and Gather (1993) described outliers "as those observations which have a different distribution from some assumed distribution for the non-outliers." Fieller (1993) notes that "there are two common themes to proposed definitions: (a) outliers are extreme observations in the sample; and (b) they are observations that are sufficiently extreme as to have an apparently low probability of occurrence or are surprising in some other way, even when adjudged as the extremes of the sample". Seaver and Triantis (1995) demonstrated that "outliers and leverage points not only affect the magnitude of the derived technical efficiency measures but also represent efficient or inefficient production performance once measurement errors have been removed from the data set". Grosskopf (1996) introduced that "outliers due to measurement error will have their most extreme effect if measurement error falsely results in observations determining the frontier, since that will then affect the efficiency scores of all the observations for which that observation forms the basis, causing their inefficiency to be overstated". Fox et al. (2004) defined that "outliers are observations which are different, in some sense, from the other observations in the sample". They also described an observation as a "scale outlier" if it is relatively larger or smaller in all, or has many dimensions of other observations, and called a "mix outlier" if it has an unusual combination in terms of the size of vector elements relative to other firms.

In this paper, the presence of outliers is classified in four (4) types. The first type is the outcome of recording/measurement errors. This type is commonly used to introduce outliers in the DEA literature.

The second type is from the lack of selecting inputs and/or outputs. In other words, if inputs/outputs of DMUs are not correctly selected, it may cause presenting outliers even if there is no measurement/recording error in data. In this type, adding an appropriate input/output may rectify the issue. Note that, drawbacks of selecting inputs and outputs are not from the structure of DEA, however, DEA is able to measure the impact of each factor on the measured efficiency score by applying KAM.

The third one is when there is a non-homogenous DMU among observations. If data are exact and the inputs/outputs are appropriately selected, a non-homogenous DMU may clearly cause outliers. For example, every airport should have the standard infrastructures according to the International Civil Aviation Organization (ICAO) documents; however, an airport in a historical tourist city is not appropriately homogenous with a military airport in a small village which may have been used for commercial transportation.

The fourth type is when a DMU, in comparison with another DMU, uses a much greater (bit less) value of input and only produces a bit greater (much less) value of output. For instance, DMU $B(10; 7.1)$ in comparison with DMU $A(2; 7)$ uses five times of input values to produce one tenth more output values only. We call this type of outliers as 'Near and Far Data' (NFD), because, for instance, the outputs of $A$ and $B$ are too close, but their inputs are too far. Therefore, if the first three types are exactly suitable, there may still be NFDs to effect on the measurement.

If all of the above four types are suitable, there is no appropriate reason to call a DMU an outlier even if the DMU has very extreme values. In this case, the relative efficiency meaning in DEA is quite enough and logical to rank and benchmark DMUs.

The above four types of outliers presence can be summarized as follows:

I. Existing measurement/recording errors.
II. Existing inappropriate/missing inputs/outputs.
III. Existing non-homogenous DMUs.
IV. Existing NFDs.

Now, let us apply 0-KAM and $\epsilon$-KAM ($\epsilon \in \mathbb{R}_+^{m+p}$). The following definition to identify outliers is proposed as follows:

**Definition 3**: A DMU is called outlier in a sample when at least one of the following cases is satisfied:

i. The technical efficiency score of DMU with 0-DF is much greater than most of DMUs' technical efficiency scores.
ii. The best technical efficiency score of DMU with $\epsilon$-DF is much greater than most of DMUs' best technical efficiency scores.
iii. The best technical efficiency score of DMU with $\epsilon$-DF moderately decreases in comparison with its technical efficiency score.
iv. It is a technically efficient DMU, and has a great sensitivity score with $\epsilon$-DF.

Note that, it is possible to use descriptive statistics such as 'deciles' and/or statistical dispersion such as 'standard deviation' instead of using the phrase 'much greater' in

Definition 3, however, as it is illustrated in the next section, we suggest depicting a diagram inclusive of polygon charts of technical efficiency, best technical efficiency, lowest efficiency and sensitivity indexes which sorted by technical efficiency index, and then reading the results from the slopes of line segments easily.

## 4. The methodology of detecting outliers

Let us suppose 100 DMUs depicted in Figure 1 with one input and one output when at least one of the types (I)-(IV) happen. This example was selected to be similar to the work of Simar (2007) and Simar and Zelenyuk (2011) which claimed the weakness of DEA to detect outliers. From the figure, it is clear that there are some outliers and extreme values which effect on the frontier.

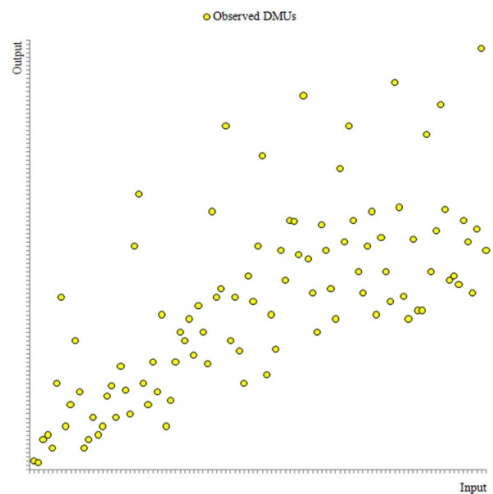

Figure 1: A sample of 100 DMUs.

Now, let us apply $\epsilon$-KAM when $w_1^-$ and $w_1^+$ are $1/x_{l1}$ and $1/y_{l1}$, $\varepsilon = 0.1$ (in order to have a clear figure) and $\varepsilon_1^-$ and $\varepsilon_1^+$ are $0.1 \times x_{l1}$ and $0.1 \times y_{l1}$, respectively. The $\epsilon$-KAM and $\epsilon$-DF in this case are denoted by 0.1-KAM and 0.1-DF, respectively.

Figure 2 depicts the polygon charts of technical efficiency scores with 0-DF and 0.1-DF, arranged from the most 0-KAM scores to the least.

The first 18 DMUs in Figure 2 were labeled A-R. DMUs A-F were technically efficient and other DMUs were inefficient. There was also a dramatic fall in the technical efficiency score with 0-DF form F to R, whereas it slightly decreased for other DMUs. Moreover, only A was KAM efficient with 0.1-DF by Definition 2 and other DMUs were inefficient with 0.1-DF. There was also a sharp decline in the best technical efficiency score with 0.1-DF from A-K. The best efficiency scores with 0.1-DF for these 18 DMUs except A, L, N, O and Q were sharply less than their corresponding technical efficiency scores with 0-DF. From these illustrations and the parts (i), (ii) and (iii) of Definition 3, KAM suggests A-R as outliers. Indeed, they may be from measurement errors, non-homogeneity, and bad selecting of factors (if there are any).

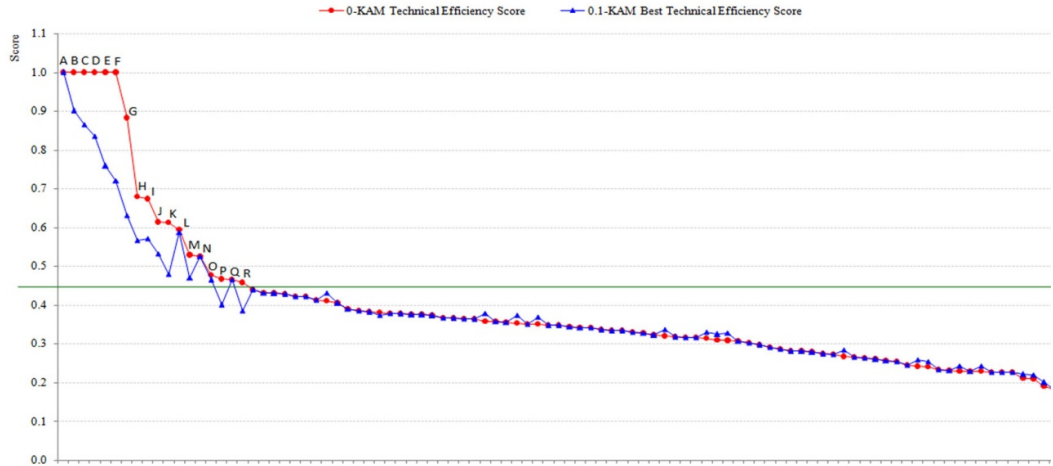

Figure 2: KAM outcomes for 100 DMUs.

Now let us see which DMUs are A-R in the sample. Figure 3 depicts DMUs A-R among observed DMUs.

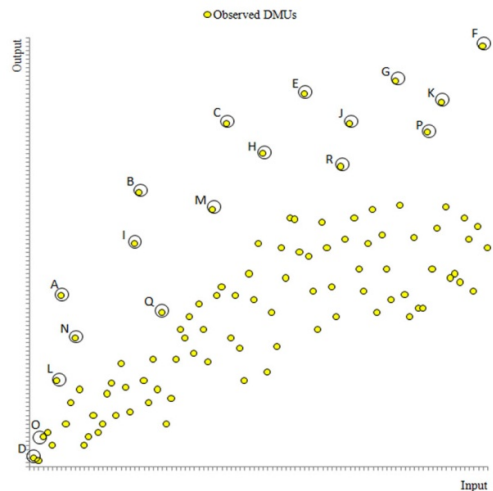

Figure 3: Detecting outliers by KAM.

As can be seen, KAM nicely detects outliers while it simultaneously ranks and benchmarks them without the computational complexities of current hybrid technologies. We also selected the half values of minimum input and minimum output of DMUs to define $\varepsilon_1^-$ and $\varepsilon_1^+$, respectively, in order to have the same commensurate effects on the Farrell frontier for evaluating each DMU (Khezrimotlagh *et al.* 2013a), and again found DMUs A-R as outliers according to the parts (i), (ii) and (iii).

The above finding represents the validity of KAM, and rejects the mentioned claim of DEA's weakness to detect outliers in Simar (2007) and Simar and Zelenyuk (2011). Moreover, the methodology is quite simple, which is one of the robust advantages of KAM in comparison with current methodologies.

Now, let us exclude DMUs A-R from the sample, and estimate the production frontier of remainder DMUs by 0.1-KAM. Figure 4 depicts the outcomes of 0-KAM and 0.1-KAM. We again labeled the first 26 DMUs by 'a' - 'z' as shown in the figure.

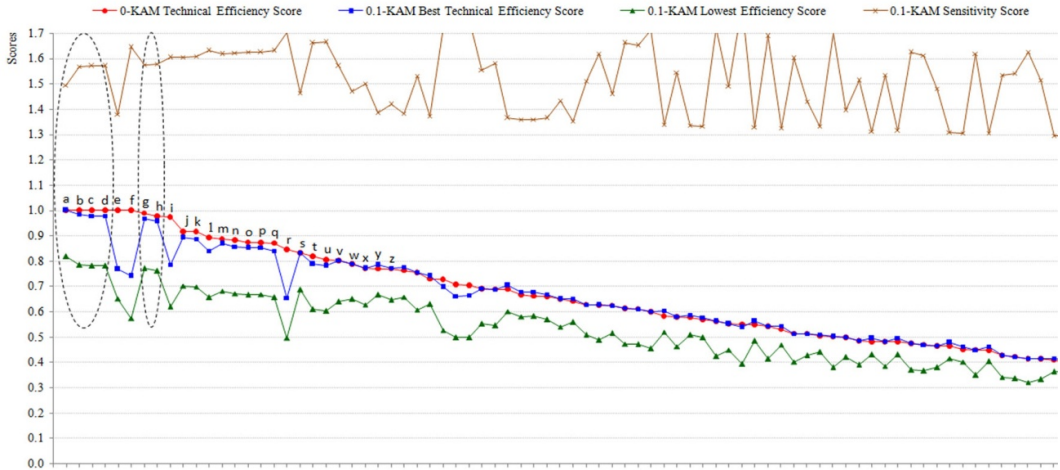
Figure 4: Detecting best efficient DMUs by KAM.

The polygon chart of technical efficiency score with 0-DF is gradually decreased. As it depicted with two dotted vertical ellipses, DMUs 'a'-'d', 'g' and 'h' have simultaneously a good value of technical efficiency score with 0-DF, a good value of best technical and lowest efficiency score with 0.1-DF, a lesser difference between technical efficiency score with 0-DF and 0.1-DF, and not a great value of sensitivity score with 0.1-DF among other DMUs. Therefore, 0.1-KAM suggests these DMUs as the best performers among remainder DMUs.

0.1-KAM also suggests DMUs 'e', 'f', 'i' and 'r' as outliers, because there is a rapid difference between their technical efficiency scores with 0-DF and 0.1-DF.

Now, let us depict where DMUs 'a'-'z' are in the sample by Figure 5. KAM properly arranges DMUs, and identifies the best performers which are DMUs 'a' - 'd', 'g' and 'h', while it detects NFDs, that is, DMUs 'e', 'f', 'r' and 'r'.

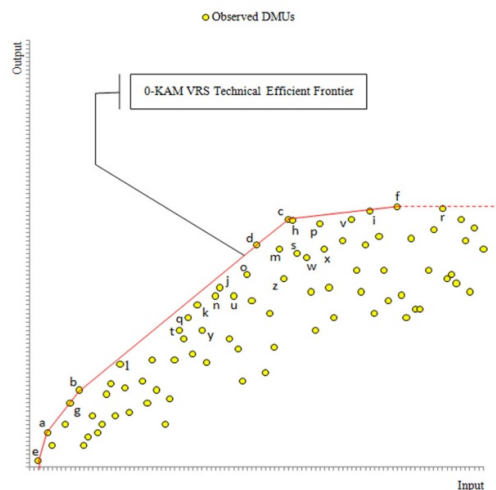
Figure 5: Detecting NFDs by KAM.

KAM also benchmarks DMUs appropriately. Figure 6 depicts 0.1-KAM technical and lowest targets. The technical efficient targets with 0.1-DF shown by small triangles are on the Farrell frontier. As can be seen, 0.1-KAM benchmarks NFDs 'f', 'i', and 'r' to DMU 'c'.

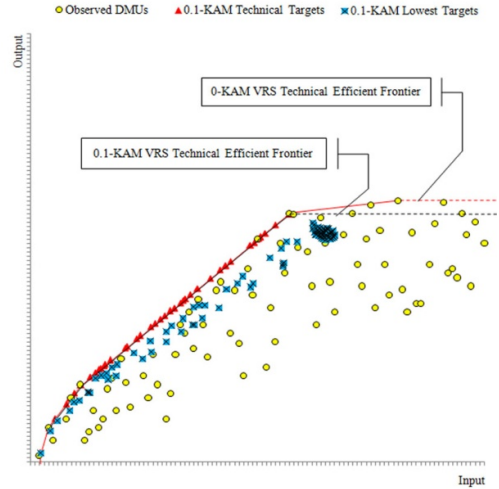
Figure 6: KAM benchmarking with 0.1-DF.

The 0.1-KAM lowest targets, shown by small multiplication signs, illustrate the effects on the Farrell frontier by introducing 10 percentage errors in DMUs' input and output values. Indeed, the effects on the Farrell frontier were commensurately considered for each evaluated DMU by selecting $\varepsilon_1^-$ and $\varepsilon_1^+$ as **0.1 × $x_{l1}$** and **0.1 × $y_{l1}$**, respectively, and can be considered as variations of the sample, which was one of the concerns of Grosskopf (1996).

We also applied KAM when $\varepsilon_1^-$ and $\varepsilon_1^+$ were **0.5 × min{$x_{i1}$: $x_{i1} \neq 0, i = 1,2,…,100$}** and **0.5 × min{$y_{i1}$: $y_{i1} \neq 0, i = 1,2,…,100$}**, respectively, and $w_1^- = 1/x_{l1}$ and $w_1^+ = 1/y_{l1}$. In this case, the effect on the Farrell frontier was the same for evaluated DMUs. There were no significant differences to arrange and benchmark DMUs. DMU 'a' still had the first rank, but it was suggested as an outlier due to have a very great sensitivity score by 0.5-DF. Moreover, the rank of technical efficient DMUs were respectively as 'a', 'c', 'd', 'b', 'f' and 'e'.

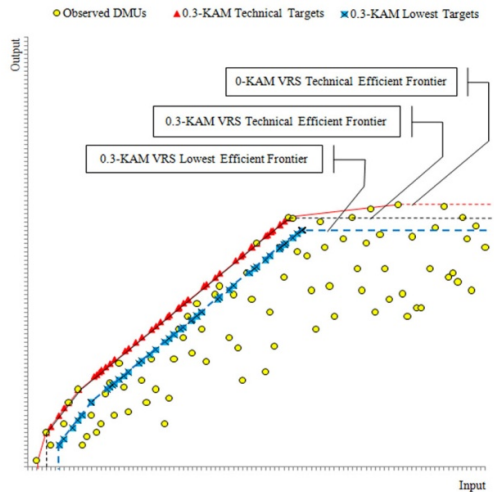
Figure 7: KAM benchmarking with 0.3-DF.

In order to have a clear figure, $\varepsilon_1^-$ and $\varepsilon_1^+$ were defined the same as 0.3 of unit and KAM was applied. Figure 7 depicts the 0.3-KAM benchmarks. From the figure, 0.3-KAM technical targets are between DMUs 'a' to 'c' on the Farrell frontier, which is the economical part of the Farrell frontier with 0.3-DF. There is no NFD in 0.3-KAM targets, too.

The methodology was also changed and supposed that the DMUs in Figure 5 are considered with one of the following assumptions: (1) DMU N from Figure 1 is added. (2) DMU M from Figure 1 is added. (3) DMU P from Figure 1 is added. For the first assumption, KAM identified DMU N as an outlier by case (ii) of Definition 3. M was known as an outlier by case (i) of Definition 3 for the second assumption, and for the last assumption, case (iv) of Definition 3 suggests P as an outlier.

These outcomes illustrate the advantages of KAM to detect extreme values without the computational complexity of the current hybrid methodologies. The above numerical example had only one input and one output, however, applying KAM to detect outliers does not have limitations for multiple inputs and multiple outputs (See Khezrimotlagh *et al.* 2012b).

## 5. Conclusion

This paper clearly illustrates how DEA is easily able to detect outliers by the technique of KAM without any additional procedures. The paper also classifies the presence of outliers into four (4) types, and proposes a definition to identify outliers in DEA.

## References

1. Andersen, P., Petersen, N.C., 1993. A procedure for ranking efficient units in data envelopment analysis. Management Science, 39 (10) 1261-1264.
2. Andrews, D.F., Pregibon, D., 1978. Finding the Outliers that Matter, Journal of the Royal Statistical Society, Series B (Methodological), 40 (1) 85-93.
3. Banker RD (1996). Hypothesis tests using data envelopment analysis. Journal of productivity analysis, 7 (2-3) 139-159.
4. Banker, R.D., Chang, H., 2006. The super-efficiency procedure for outlier identification, not for ranking efficient units. European Journal of Operational Research, 175 (2) 1311-1320.
5. Barnett, V., Lewis, T., 1984. Outliers in statistical data. Volume 3, Wiley, New York.
6. Cazals, C., Florensa, J.P., Simar, L., 2002. Nonparametric frontier estimation: a robust approach. Journal of Econometrics, 106 (1) 1-25.
7. Charnes, A., Cooper, W.W., Rhodes, E., 1978. Measuring the inefficiency of decision-making units. European Journal of Operational Research, 2 (6) 429-444.
8. Charnes, A., Neralic, L., 1990. Sensitivity Analysis of the Additive Model in Data Envelopment Analysis. European Journal of Operational Research, 48 (3) 332-341.
9. Chen, W.C., Johnson, A., 2000. Detecting efficient and inefficient outliers in data envelopment analysis. Available at SSRN 929971.
10. Chen, W.C., Johnson, A.L., 2010. A unified model for detecting efficient and inefficient outliers in data envelopment analysis. Computers \& Operations Research, 37 (2) 417-425.
11. Chung, W., (2011). Review of building energy-use performance benchmarking methodologies. Applied Energy, 88 (5) 1470-1479.
12. Coelli, T.J., Rao, D.S.P., O'Donnell, C.J., Battese, G.E., 2005. An Introduction to Efficiency and Productivity Analysis. Second Edition, Springer, New York.
13. Cooper, W.W., Seiford, L.M., Tone, K., 2007. Data Envelopment Analysis, A comprehensive Text with Models, Applications, References and DEA-Solver Software. 2nd Edition, Springer, New York.
14. De Jorge Moreno, J., Sanz-Triguero, M., 2011. Estimating technical efficiency and bootstrapping Malmquist indices: Analysis of Spanish retail sector. International Journal of Retail & Distribution Management, 39 (4) 272-288.
15. De Sousa, M.D.C.S., Stosic, B., 2005. Technical efficiency of the Brazilian municipalities: correcting nonparametric frontier measurements for outliers. Journal of Productivity analysis, 24 (2) 157-181.

16. De Witte, K., Marques, R.C., 2010. Influential observations in frontier models, a robust non-oriented approach to the water sector.~Annals of Operations Research, 181 (1) 377-392.
17. Drake, L., Simper, R., 2002. X-efficiency and scale economies in policing: a comparative study using the distribution free approach and DEA. Applied Economics, 34 (15) 1859-1870.
18. Dula, J.H., 2008. A computational study of DEA with massive data sets. Computers \& Operations Research, 35 (4) 1191-1203.
19. Emrouznejad, A., De Witte, K., 2010. COOPER-framework: A unified process for non-parametric projects. European Journal of Operational Research, 207 (3) 1573-1586.
20. Estache, A., De La Fe, B.T., Trujillo, L., 2004. Sources of efficiency gains in port reform: a DEA decomposition of a Malmquist TFP index for Mexico. Utilities policy, 12 (4) 221-230.
21. Fieller, N.R.J., 1993. Discussion of Davies and Gather (1993). Journal of the American Statistical Association, 88 794-795.
22. Fried, H.O., Lovell, C.K., Schmidt, S.S., Yaisawarng, S., 2002. Accounting for environmental effects and statistical noise in data envelopment analysis. Journal of productivity Analysis, 17 (1-2) 157-174.
23. Fox, K.J., Hill, R.J., Diewert, W.E., 2004. Identifying outliers in multi-output models. Journal of Productivity Analysis, 22 (1-2) 73-94.
24. Grosskopf, S., Valdmanis, V., 1987. Measuring hospital performance: A non-parametric approach. Journal of Health Economics, 6 (2) 89-107.
25. Grosskopf, S., 1996. Statistical inference and nonparametric efficiency: A selective survey. Journal of Productivity Analysis, 7 (2-3) 161-176.
26. Gunst, R.F., Mason, R.L., 1980. Regression analysis and its application: a data-oriented approach, Volume 2, CRC Press, New York.
27. Davies, L., Gather, U., 1993. The identification of multiple outliers. Journal of the American Statistical Association, 88 (423) 782-792.
28. Dusansky, R., Wilson, P.W., 1995. On the relative efficiency of alternative modes of producing a public sector output: the case of the developmentally disabled. European Journal of Operational Research, 80 (3) 608-618.
29. Hodge, V.J., Austin, J., 2004. A survey of outlier detection methodologies. Artificial Intelligence Review, 22 (2) 85-126.
30. Horta, I.M., Camanho, A.S., Moreira da Costa, J., 2012. Performance assessment of construction companies: A study of factors promoting financial soundness and innovation in the industry. International Journal of Production Economics, 137 (1) 84-93.
31. Johnson, A.L., McGinnis, L.F., 2008. Outlier detection in two-stage semiparametric DEA models. European Journal of Operational Research, 187 (2) 629-635.
32. Khezrimotlagh, D., Salleh S. and Mohsenpour, Z. 2012a, A New Method in Data Envelopment Analysis to Find Efficient Decision Making Units and Rank Both Technical Efficient and Inefficient DMUs Together, Applied Mathematical Sciences, 6(93):4609-4615.
33. Khezrimotlagh, D., Mohsenpour, Z. and Salleh, S. 2012b, Cost-Efficiency by Arash Method in DEA, Applied Mathematical Sciences, 6(104):5179-5184.
34. Khezrimotlagh, D., Salleh, S. and Mohsenpour, Z. 2012c, Arash Method and Uncontrollable Data in DEA, Applied Mathematical Sciences, 6(116):5763-5768.
35. Khezrimotlagh, D., Mohsenpour, Z. and Salleh, S. 2012d, Comparing Arash Model with SBM in DEA, Applied Mathematical Sciences, 6(104):5185-5190.
36. Khezrimotlagh, D., Salleh S. and Mohsenpour, Z. 2012e, Airport Efficiency with Arash Method in DEA, Journal of Basic and Applied Scientific Research, 2(12):12502-12507.
37. Khezrimotlagh, D., Salleh, S. and Mohsenpour, Z. 2012f, Benchmarking Inefficient Decision Making Units in DEA, Journal of Basic and Applied Scientific Research, 2(12):12056-12065.
38. Khezrimotlagh, D., Salleh, S., Mohsenpour, Z., 2013a. A New Method for Evaluating Decision Making Units in DEA. Journal of the Operational Research Society, 65 (1) 694–707.
39. Khezrimotlagh, D., Salleh, S., Mohsenpour, Z., 2013b. A New Robust Mixed Integer Valued Model in DEA. Applied Mathematical Modelling, 37 (24) 9885-9897.
40. Khezrimotlagh, D., Salleh, S. and Mohsenpour, Z. 2013c, Nonlinear Arash Model in DEA, Research Journal of Applied Sciences, Engineering and Technology, 5(17): 4268-4273.
41. Khezrimotlagh, D., Mohsenpour, P., Salleh, S. and Mohsenpour, Z. 2013d. A Review on Arash Method in Data Envelopment Analysis. In: Banker B., A. Emrouznejad, H. Bal, I. Alp, M. Ali Gengiz (2013), Data Envelopment Analysis and Performance Measurement: Proceedings of the 11th International Conference of DEA, June 2013, Samsun, Turkey, ISBN: 9781854494771.
42. Khezrimotlagh, D. 2014a, Profit Efficiency by Kourosh and Arash Model, Applied Mathematical Sciences, 8(24):1165-1170.

43. Khezrimotlagh, D., 2014b. How to Select an Epsilon in Kourosh and Arash Model, Emrouznejad, A., R. Banker R., S. M. Doraisamy and B. Arabi (2014), Recent Developments in Data Envelopment Analysis and its Applications, Proceedings of the 12th International Conference of DEA, April 2014, Kuala Lumpur, Malaysia, ISBN:9781854494870, pp.230-235.
44. Kuosmanen, T., Johnson, A.L., 2010. Data envelopment analysis as nonparametric least-squares regression. Operations Research, 58 (1) 149-160.
45. Mahlberg, B., Raveh, A., 2012. Co-plot: a useful tool to detect outliers in DEA. Available at SSRN 1999370.
46. Mohsenpour, P., Munisamy, S. and Khezrimotlagh, D. (2013), Revenue Efficiency and Kourosh Method in DEA, Applied Mathematical Sciences, 7(140):6961-6966.
47. Ondrich, J., Ruggiero, J., 2002. Outlier detection in data envelopment analysis: an analysis of Jackknifing. Journal of the Operational Research Society, 53 (3) 342-346.
48. Pastor, J.T., Ruiz, J.L., Sirvent, I., 1999. A statistical test for detecting influential observations in DEA. European Journal of Operational Research, 115 (3) 542-554.
49. Ray, S.C., 2004. Data Envelopment Analysis: Theory and Techniques for Economics and Operations Research. Cambridge University Press, New York.
50. Seaver, B.L., Triantis, K.P., 1995. The Impact of Outliers and Leverage Points for Technical Efficiency Measurement Using High Breakdown Procedures. Management Science, 41 (6) 937-956.
51. Seiford, L.M., Thrall, R.M., 1990. Recent developments in DEA: the mathematical programming approach to frontier analysis. Journal of econometrics, 46 (1) 7-38.
52. Sengupta, N., 1991. Managing common property: irrigation in India and the Philippines. Sage Publications India Pvt Ltd, New Delhi.
53. Sexton, T.R., 1986. The methodology of data envelopment analysis. New Directions for Program Evaluation 1986 (32) 7-29.
54. Simar, L., 1996. Aspects of statistical analysis in DEA-type frontier models. Journal of Productivity Analysis, 7 (2-3) 177-185.
55. Simar, L., Wilson, P.W., 2000. A general methodology for bootstrapping in non-parametric frontier models. Journal of applied statistics, 27 (6) 779-802.
56. Simar, L., 2003. Detecting outliers in frontier models: a simple approach. Journal of Productivity Analysis 20 (3) 391-424.
57. Simar, L., 2007. How to improve the performances of DEA/FDH estimators in the presence of noise? Journal of Productivity Analysis, 28 (3) 183-201.
58. Simar, L., Wilson, P.W., 2008. Statistical inference in nonparametric frontier models: recent developments and perspectives. The measurement of productive efficiency and productivity growth 421-521. At Fried HO, Lovell CK and Schmidt SS (2008). The measurement of productive efficiency and productivity growth, Oxford University Press, New York.
59. Simar, L., Zelenyuk, V., 2011. Stochastic FDH/DEA estimators for frontier analysis. Journal of Productivity Analysis, 36 (1) 1-20.
60. Tran, N.A., Shively, G., Preckel, P., 2010. A new method for detecting outliers in Data Envelopment Analysis. Applied Economics Letters, 17 (4) 313-316.
61. Wilson, P.W., 1993. Detecting outliers in deterministic nonparametric frontier models with multiple outputs. Journal of Business and Economic Statistics, 11 (3) 319-323.
62. Wilson, P.W., 1995. Detecting influential observations in data envelopment analysis. Journal of Productivity Analysis, 6 (1) 27-45.